\documentclass{article}

\usepackage[english,russian]{babel}

\usepackage{amsmath}
\usepackage{amsfonts}
\usepackage{amssymb,latexsym}
\usepackage{amsthm}

\newtheorem{theorem}{Теорема}[section]
\newtheorem{definition}{Определение}[section]
\newtheorem{lemma}{Лемма} [section]

\newtheorem{corollary}{Следствие}

\newtheorem{{thebibliography}}{REf}

\begin{document}

\author{М.З. Двейрин, А.С. Левадная}

\title{Обобщенный порядок и обобщенный тип целой функции в терминах ее наилучших приближений}

\date{}



\maketitle

Аннотация.

В статье выясняется связь между обобщенным порядком и обобщенным
типом целой функции бесконечного порядка и скоростью ее наилучшей
полиномиальной аппроксимации для большого семейства банаховых
пространств функций, аналитических в единичном круге. Найдены
соотношения, определяющие обобщенные порядок и тип целой функции
через последовательность  ее наилучших приближений. Полученные
результаты являются обобщением более ранних результатов Редди, Д.
Сато, И.И. Ибрагимова и Н.И. Шихалиева, С.Б. Вакарчука, Р. Мамадова.

Annotation.

 The paper explores connection between the
generalized order and the generalized type of an entire function and
the speed of the best polynomial approximation in the unit disk. The
relations which define the generalized order and the generalized
type of an entire function through the sequence of its best
approximations, have been found. The results were obtained by
generalization previous results of A.R. Reddy, D. Sato, I.I.
Ibragimov and N. I. Shyhaliev, S. B. Vakarchyk, R. Mamadov.


Ключевые слова и фразы. Целая функция, наилучшее приближение,
обобщенный порядок целой функции, обобщенный тип целой функции.

\section{Введение}

В данной работе в качестве $X$ рассматривается линейное
нормированное пространство, образованное аналитическими в единичном
круге $D$ функциями, имеющими конечную норму
$\parallel\cdot\parallel$. При этом будем предполагать, что
$\parallel\cdot\parallel$ помимо обычных свойств нормы удовлетворяет
также условиям:
\begin{equation}
\label{trivial} i) \quad \, \|f(\cdot e^{it})\| \, \equiv
\|f(\cdot)\| \,
\end{equation}
для всех $t \in \mathbb{R}$ и $f \in X;$

\begin{equation}
ii) \quad \, \|f(\cdot)\| \,<\infty  \,
\end{equation}
для   любой целой функции (т.е. пространство $X$ содержит все целые
функции);

\begin{equation}
iii) \quad \, \|  \frac
{1}{2\pi}\int\limits_{0}^{2\pi}f(ze^{it})g(t) \,dt \| \leq \, \frac
{1}{2\pi}\int\limits_{0}^{2\pi}|g(t)| \,dt \,\, \|f(\cdot)\|
\end{equation}
для   любых функций $f \in X$  и  $g \in L[0; 2\pi]$   (иначе
говоря, для любых $f \in X$  и  $g \in L[0; 2\pi]$ \, $\|f \ast g \|
\leq \| f \|\, \|g\|_{L[0; 2\pi]}$).

Этим требованиям удовлетворяет норма в целом ряде функциональных
пространств, являющихся объектом многочисленных исследований (авторам
неизвестны примеры пространств, в которых выполняются условия i),
ii) и при этом не выполняется условие iii) ). Приведем некоторые из
них.

 1) Пространство $B$ функций, аналитических в единичном круге $\mathbb{D}$ и
непрерывных на его замыкании $\overline {\mathbb{D} }$ с нормой

\[
\|f\|=\max\limits_{z\in \overline{\mathbb{D}}}|f(z)|<\infty \,.
\]

 2) Пространства Харди $H_p$ ($p \geq 1$) функций, аналитических в круге
$\mathbb{D}$   с нормой

\[
\|f\|=\sup \limits_{0<r<1}M_{p} \,(f,r),\quad M_{p} \,(f,r):=
\left(\frac {1}{2\pi}\int\limits_{0}^{2\pi}|f(re^{it})|^p
\,dt\right)^\frac{1}{p}\;,\quad p\in [1; \infty);
\]

\[
\|f\|=\sup\limits_{z\in \mathbb{D}}|f(z)|\,,\, \quad p=\infty.
\]

 3) Пространства Бергмана $H_p^\prime$ функций, аналитических в круге
$\mathbb{D}$ при $p\in [1; \infty)$ с нормой

\[
\|f\|=\left(\frac {1}{\pi} \int \limits_{\;\; z \in D}\int
|f(x+iy)|^p \;dxdy\right)^\frac{1}{p}\;,
\]
и обобщенные (весовые) пространства Бергмана $H_{p, \, \rho}^\prime$
функций, аналитических в круге $\mathbb{D}$ при $p\in [1; \infty)$ с
нормой
\[
\|f\|=\left(\frac {1}{\pi} \int \limits_{\;\; z \in D}\int
|f(x+iy)|^p \;\rho(|z|)dxdy\right)^\frac{1}{p}\;
\]
и радиальным весом $\rho(|z|)$.

 4) Пространства  $A_p, \quad p \in (0;1)\; $ функций,
аналитических в круге $\mathbb{D}$  с нормой

\[
\|f\|=\int\limits_{0}^{1}(1-r)^{\frac{1}{p}-2}M_{1} \,(f,r) \,dr\, ,
\]
впервые изучавшиеся Харди и Литллвудом \cite{HL2} и позднее
Ромбергом, Дюреном и Шилдсом  \cite{DRS}.

 5) Пространства  $\mathcal{B}_{p,\,q,\,\lambda}\, , \quad
0<p<q\leq\infty, \quad \lambda > 0,\; $ функций, аналитических в
круге $\mathbb{D}$ с нормой

\[
\|f\|=\left\{ \int \limits_{0}^{1}(1-r)^{\lambda
 \, p \;q \,(q-p)^{-1}}M_q^\lambda(f,r) \,dr\right\}^\frac{1}{\lambda},\, \quad
\lambda<\infty,
\]

\[
\|f\|=\sup\limits_{0<r<1}\left\{ (1-r)^{
 \, p \;q \,(q-p)^{-1}}M_{q}\,(f,r)\right\} \, ,\quad \lambda=\infty,
\]
введенные {Харди и Литллвудом в работе \cite{HL2} (см. также
\cite{Gvar}).

 6) Пространства со смешанной нормой $H^{p,\,q,\,\alpha} , \,
(p, q \geq 1, \quad\,  \alpha >0),$ образованные функциями,
аналитическими в круге $\mathbb{D}$ с конечной нормой

\[
\|f\|=\left\{ \int \limits_{0}^{1}(1-r)^{q\alpha -1}M_p^q(f,r)
\,dr\right\}^\frac{1}{q} \, ,    \quad q<\infty,
\]

\[
\|f\|=\sup\limits_{0<r<1}\left\{ (1-r)^{
 \alpha}M_{p}\,(f,r)\right\} \, ,\quad q=\infty,
\]
введенные Харди и Литллвудом в работе  \cite{HL2}.    Заметим, что
пространства   $H^{p,\,q,\,\alpha}$  и
$\mathcal{B}_{p,\,q,\,\lambda}\,$  отличаются лишь способом введения
параметров.

 7) Пространство $\quad BMOA$ \cite{Swed}, состоящее из функций $f\in H_{1}$ с нормой
\[
\|f\|=\sup\limits_{I}\int \limits_{I}
|f(\zeta)-f_{I}|d\sigma(\zeta)\,,
\]
где $f(\zeta)$ - граничные значения функции $f(z)$ на единичной
окружности, а $f_{I}$ - среднее арифметическое значение функции
$f(\zeta)$ на дуге $I$.

 8) Пространства типа Блоха $\mathcal{B}_{\alpha}\,$,  $\alpha\in(0, \infty)$, состоящие из функций, аналитических в  $D$ с конечной нормой
\[
\|f\|=|f(0)|+ \sup\limits_{z\in D} (1-|z|^{2})^{\alpha}
|f^{\prime}(z)|\,.
\]
Пространства  $\mathcal{B}_{\alpha}\,$  являются банаховыми
\cite{Zhu},  при $\alpha=1$   $\mathcal{B}_{\alpha}$ совпадает с
пространством Блоха $\mathcal{B}$.

 9) Введенные Е.М. Дынькиным \cite{Dink} пространства  $\mathcal{A}_{p,q}^s (\mathbb{D})\,$
 функций, являющиеся аналогами классов О.В. Бесова $\mathcal{B}_{p,q}^s
[-1;\;1]\,$.  Эти пространства образованы функциями $f\in H_{p}$,
$p\in [1; \infty]$  с нормой

\[
\|f\|=\left\{ \int \limits_{0}^{1}\left (  \frac{\omega_{m}(f, \,
t)_{p}}{t^{s}} \right)^{q}\frac{dt}{t}\right\} ^\frac{1}{q} \, +
\sup \limits_{0<r<1}M_{p} \,(f,r).
\]

Здесь  $q\in [1; \infty], \, s>0, \, m>s$ - натуральное число,
$\omega_{m}(f, \, t)_{p}$ - $m$-ый модуль гладкости в пространстве
$L_p\,$ функции  $f(e^{i\cdot})$, представляющей собой радиальные
предельные значения $f$.  Случай $q=\infty$ трактуется традиционно.

10) Обобщенные пространства  Дирихле $\mathcal{D}_{p}(\alpha)\,$
функций, аналитических в  $\mathbb{D}$,  с нормой
\[
\|f(z)\|=\left(
\sum\limits_{k=0}^{\infty}|c_{k}|^{p}\,\alpha_{k}\right)^{1/p}\,,
\]
где   $ c_{k}=c_{k}(f)$ - коэффициенты Тейлора функции $f$,  $p\geq
1$,  ${\alpha}=\{\alpha_{k}\}$ - фиксированная последовательность
положительных чисел с условиями
$$ \limsup \limits_{k\rightarrow\infty}\left(\alpha_{k}\right)^{\frac{1}{k}}< \infty, \quad
 \liminf \limits_{k\rightarrow\infty}\left(\alpha_{k}\right)^{\frac{1}{k}}\geq 1 . $$

Отметим, что приведенные выше примеры функциональных  пространств со
свойствами i), ii) и iii) не исчерпывают их многообразия.

Напомним общепринятые определения основных характеристик целой
функции. В дальнейшем будем использовать функцию $\ln\ln..\ln x$,
где $x$ логарифмируется $q$ раз. Введем для нее обозначение:
$\ln^{(q)}x:=\ln\ln..\ln x$, $\ln^{(0)}x:=x$. Будем также обозначать
$$M(f,r):=\sup_{0<|z|<r}|f(z)|.$$

\begin{definition}\label{por}
Порядок роста целой функции равен
$$\rho:=\limsup_{r\rightarrow\infty}{\frac{\ln\ln M(f,r)}{\ln r}}.$$
\end{definition}

\begin{definition}\label{type}
Тип целой функции равен (если $0<\rho<\infty$) $$\sigma:=
\limsup_{r\rightarrow\infty}{\frac{\ln M(f,r)}{ r^\rho}}.$$
\end{definition}

\begin{definition}\label{por1}
Обобщенный порядок роста $\rho_q$ индекса q (q-порядок) целой
функции равен (в случае, когда $\rho=\infty$)
$$\rho_q:=\limsup_{r\rightarrow\infty}{\frac{\ln^{(q)} M(f,r)}{\ln r}}$$
\end{definition}
где $q$ - натуральное число, удовлетворяющее условиям
$\rho_{q-1}=\infty$, $\rho_q<\infty$.

\begin{definition}\label{type1}
Обобщенный тип $\sigma_q$ индекса q (q-порядок) целой функции равен
(для $0<\rho_q<\infty$)
$$\sigma_{q}:= \limsup_{r\rightarrow\infty}{\frac{\ln^{(q-1)}
M(f,r)}{ r^{\rho_q}}}.$$
\end{definition}

Данные определения обобщенного порядка и обобщенного типа введены
Сато Д. в \cite{Sato} и Редди А.Р. в \cite{Rd1} (в случае $q=2$ мы
будем использовать обозначения $\rho$ вместо $\rho_2$ и $\sigma$
вместо $\sigma_2$). Введение q-порядка и q-типа позволяет различать
скорость роста целой функции в случае, когда $\rho=\infty.$

В частности, в статье Сато Д. \cite{Sato}  были получены формулы,
которые связывают величины $\rho_q$ и $\sigma_q$ с тейлоровскими
коэффициентами $c_n$ целой трансцендентной функции $f$:

$$\rho_q=\limsup_{n\rightarrow\infty}\frac{n\ln^{(q-1)}
n}{-\ln|c_n|}$$ если $q=2,3,...$ и

$$\sigma_q=\begin{cases}
\ \frac{1}{e\rho_q}\limsup_{n\rightarrow\infty}n|c_n|^\frac{\rho_q}{n} ,&\text{если $q=2$;}\\
\limsup_{n\rightarrow\infty}\ln^{(q-2)}n\cdot|c_n|^\frac{\rho_q}{n},&\text{если
$q=3,4,..$.}
\end{cases}$$

Связи между ростом максимума модуля целых функций или функций,
аналитических в круге, и наилучшим приближением изучались в работах
Редди А.Р. \cite{Rd1}, Ибрагимова И.И. и Шихалиева Н.И. \cite{Ibr1},
\cite{Ibr2}, Вакарчука С.Б. \cite{Wac}, Мамадова Р. \cite{Mamad} и
Двейрина М.З. \cite{Dmz}. Ими были получены
 соотношения, выражающие порядок и тип целой функции через
последовательность ее наилучших приближений $E_n(f)$ для
аппроксимации по различным нормам. Более полное изложение истории
исследований по данной теме можно найти в \cite{WacGir}.

Напомним, что $E_n(f)\equiv E_n(f,L_n)$ - наилучшее приближение
функции $f \in X$ элементами линейного подпространства $L_n$
определяется следующим образом:
$$E_n(f):=\inf_{p\in L_n}{\| f-p\|}$$
В качестве аппроксимирующего подпространства будем использовать
${\mathcal{P}}_n$ - совокупность алгебраических полиномов
комплексной переменной степени не выше $(n-1)$.

Приведем некоторые нужные для дальнейшего факты из статьи
 \cite{Dmz}:

\begin{theorem}\label{th1.2}
Пусть $f\in X$. Тогда условие
$$\lim_{n\rightarrow\infty}{\left(E_n(f)\right)^\frac{1}{n}}=0$$
является необходимым и достаточным для того, чтобы функция $f$ была
целой.
\end{theorem}

\begin{lemma}\label{lemma1.3}
Пусть $f \in X$ и $f(z)=\sum_{k=0}^{\infty} c_kz^k$  при $ z \in
\mathbb{D}$. Тогда $|c_n|\cdot\|z^n\|\leq E_n(f) \leq \|f\|.$
\end{lemma}

\begin{lemma}
Пусть $ f\in X \quad$ и $\, \mu_{1}:=\liminf
\limits_{n\rightarrow\infty}\,(\|z^{n}\|)^{\frac{1}{n}}$,
$\mu_{2}:=\limsup
\limits_{n\rightarrow\infty}\,(\|z^{n}\|)^{\frac{1}{n}}$. Тогда $
\mu_{1}\geq 1,  \,  \mu_{2} < \infty$.
\end{lemma}

\begin{theorem}\label{th1.4} Для того, чтобы функция $f\in X$ была целой конечного порядка
$\rho\in (0;+\infty)$ необходимо и достаточно, чтобы существовал
конечный положительный предел
\begin{equation}\label{1}
\limsup_{n\rightarrow\infty}\frac{n\ln
n}{\ln\frac{\|z^n\|}{E_n(f)}}=\alpha.
\end{equation}
При этом справедливо равенство $\alpha=\rho.$
\end{theorem}

\begin{theorem}\label{th1.5}
Пусть существует конечный предел
 $\liminf_{n\rightarrow\infty}{(\|z^n\|)^\frac{1}{n}}=\mu>0$. Для того, чтобы функция $f\in X$ была целой
конечного порядка $\rho\in (0;+\infty)$ и нормального типа
$\sigma\in (0;+\infty)$ необходимо и достаточно, чтобы
\begin{equation}\label{2}
\sigma=\frac{1}{e\rho}\limsup_{n\rightarrow\infty}n\left(\frac{E_n(f)}{\|z^n\|}\right)^\frac{\rho}{n}.
\end{equation}
\end{theorem}

В настоящей статье будут получены аналоги теорем
\ref{th1.2}-\ref{th1.5}, устанавливающие связи между наилучшими
полиномиальными приближениями функции и ее обобщенным порядком
$\rho_q$ и обобщенным типом $\sigma_q$  в случае $q>2$ (случай $q=2$
соответствует результатам, полученным ранее в \cite{Dmz}).

Всюду в дальнейшем будем обозначать

\begin{equation}\label{5.1}\mu_{1}:=\liminf
\limits_{n\rightarrow\infty}\,(\|z^{n}\|)^{\frac{1}{n}}, \,\,\,\,\,
 \mu_{2}:=\limsup
\limits_{n\rightarrow\infty}\,(\|z^{n}\|)^{\frac{1}{n}}.\end{equation}

\section{ Формулировка результатов.}

\begin{theorem}\label{th1}

Для того, чтобы функция $f\in X$ была целой обобщенного порядка
$\rho_q\in (0;+\infty)\quad (q\in \mathbb{N}, \,q\geq2) $ необходимо
и достаточно, чтобы
\begin{equation}\label{3}
\limsup_{n\rightarrow\infty}\frac{n\ln^{(q-1)}
n}{\ln\frac{\|z^n\|}{E_n(f)}}=\alpha
\end{equation}
был конечным и положительным. При этом $\alpha=\rho_q.$
\end{theorem}

\begin{theorem}\label{th2}
Пусть существует конечный предел $\lim\limits_{n\rightarrow \propto}
(\|z^{n}\|)^\frac{1}{n}=\mu $. Для того, чтобы функция $f\in X$ была
целой конечного порядка $\rho_q\in (0;+\infty)$ и типа $\sigma_q\in
(0,\infty)$ необходимо и достаточно, чтобы
\begin{equation}\label{4}
\sigma_q=\begin{cases}
\ \frac{1}{e\rho_q}\limsup_{n\rightarrow\infty}n\left(\frac{E_n(f)}{\|z^n\|}\right)^\frac{\rho_q}{n} ,&\text{если $q=2$;}\\
\limsup_{n\rightarrow\infty}\ln^{(q-2)}n\cdot\left(\frac{E_n(f)}{\|z^n\|}\right)^\frac{\rho_q}{n},&\text{если
$q=3,4,...$} \,\,.
\end{cases}
\end{equation}
\end{theorem}

 \begin{corollary}
Условие
\begin{equation}
\lim_{n\rightarrow\infty} \, \left(\frac{E_n(f)}{\|z^n\|}\right)^{\frac{1}{n}} \ln^{(q-2)}
n = 0
\end{equation}
является необходимым и достаточным для того, чтобы функция $f$  была целой некоторого обобщенного порядка  $\rho_q$  с  $\rho_q\in (0; 1)$.
\end{corollary}

Чтобы сформулировать следующее следствие нам понадобится ввести некоторые обозначения.    Пусть  $ \Omega$ - ограниченный континуум со связным дополнением в комплексной плоскости,  $0 \in \Omega$ - его внутренняя точка,  $r$ и $R$ - радиусы кругов  $D_r$ и $D_R$ с центром в точке $z=0$ и таких, что  $D_r \subset \Omega \subset D_R$;  для целой функции $f$ положим  $f_r(z):=f(rz)$  и
$$
\|f\|_\Omega := \left( \int \int \limits_{\Omega} |f(x+iy)|^{p}dxdy \right)^{\frac{1}{p}},  \quad p\geq 1
$$

($\|f\|_\Omega$ - норма функции $f$ в пространстве $E^{'}_{p}(\Omega)$, которое в случае, когда  $ \Omega$ есть замыкание области, представляет собой хорошо  известное пространство В.И. Смирнова),   $E_n(f)_\Omega$ - наилучшее приближение
функции $f$  алгебраическими полиномами  комплексной переменной степени не выше $(n-1)$ в пространстве  $E^{'}_{p}(\Omega)$.

 \begin{corollary}
Пусть функция $f \in E^{'}_{p}(\Omega), \quad p\geq 1$. Для того чтобы $f$ была целой обобщенного порядка $\rho_q\in (0; \infty)$  необходимо и достаточно, чтобы выполнялось условие
\begin{equation}
\limsup_{n\rightarrow\infty}\frac{n\ln^{(q-1)}
n}{-\ln E_n(f)_\Omega}=\rho_q.
\end{equation}
\end{corollary}

Отметим, что полученные результаты в частных случаях пространств  $\quad BMOA$, $\mathcal{B}_{\alpha}\,$,  $\mathcal{A}_{p,q}^s (\mathbb{D})\,$, $\mathcal{D}_{p}(\alpha)\,$   являются новыми.  Утверждение следствия 2 при существенных ограничениях на  $\Omega$ было получено ранее в \cite{WacGir}.

\section{ Доказательства.}

Докажем теорему \ref{th1}:

\begin{proof}
\begin{flushleft}
    \emph{Достаточность.}
\end{flushleft}

Из условия (\ref{3}) теоремы \ref{th1} следует равенство
$\lim_{n\rightarrow\infty}{\left(\frac{E_n(f)}{\|z^n\|}
\right)^\frac{1}{n}}=0$, т.е. выполнение условия теоремы
\ref{th1.2}. Cледовательно, $f$ - целая функция с $\rho=\infty$.
Обозначим ее q-порядок $\rho_q$ и положим
$\alpha_q:=\limsup_{n\rightarrow\infty}\frac{n\ln^{(q-1)}
n}{\ln\frac{\|z^n\|}{E_n(f)}}$. Тогда ввиду леммы \ref{lemma1.3}
имеем неравенство:
\begin{equation}\label{5}
\rho_q=\limsup_{n\rightarrow\infty}\frac{n\ln^{(q-1)}
n}{-\ln|c_n|}\leq \limsup_{n\rightarrow\infty}\frac{n\ln^{(q-1)}
n}{\ln\frac{\|z^n\|}{E_n(f)}}=\alpha_q.
\end{equation}

Покажем, что в условиях теоремы $\rho_q>0$.

Предположим противное, т.е. что
$\limsup_{n\rightarrow\infty}\frac{n\ln^{(q-1)} n}{-\ln|c_n|}=0$.
Тогда для произвольного положительного $ \varepsilon \in (0,\mu_1)$
найдется $ N_{\varepsilon}$ такое, что при $ n
> N_{\varepsilon}$ выполняется неравенство  $ n\ln^{(q-1)} n<-\varepsilon \ln|c_n|$. Последнее
неравенство равносильно следующему: $
|c_n|<(\ln^{(q-2)}n)^{-\frac{n}{\varepsilon}}$. Пользуясь им, оценим
$E_n(f)$ при $n > N_{\varepsilon}$. Будем считать
$N_{\varepsilon}$ столь большим, чтобы выполнялись неравенства:
$\|z^n\|\leq(\mu_2+\varepsilon)^n$,
$\|z^n\|\geq(\mu_1-\varepsilon)^n$ и
${\mu_2+\varepsilon}<{(\ln^{(q-2)}n)^{\frac1 \varepsilon}}$ при $n
> N_{\varepsilon}$, где $\mu_1$ и $\mu_2$ определяются
соотношениями (\ref{5.1}). Тогда

$$E_n(f)\leq \| \sum_{k=n}^\infty c_kz^k\| \leq
\sum_{k=n}^\infty |c_k|\cdot\|z^k\|\leq \sum_{k=n}^\infty
|c_k|\cdot(\mu_2+\varepsilon)^k\leq $$ $$\leq \sum_{k=n}^\infty
\left[\frac{\mu_2+\varepsilon}{(\ln^{(q-2)}n)^{\frac{1}{\varepsilon}}}\right]^k
\leq
\frac{1}{1-\frac{\mu_2+\varepsilon}{(\ln^{(q-2)}n)^{\frac{1}{\varepsilon}}}}\cdot\left[\frac{\mu_2+\varepsilon}{(\ln^{(q-2)}n)^{\frac{1}{\varepsilon}}}\right]^n.$$
Следовательно,
$$\frac{\|z^n\|}{E_n(f)}\geq (\mu_1-\varepsilon)^n\cdot\left[{1-\frac{\mu_2+\varepsilon}{(\ln^{(q-2)}n)^{\frac{1}{\varepsilon}}}}\right]\cdot {\left[\frac{\mu_2+\varepsilon}{(\ln^{(q-2)}n)^{\frac{1}{\varepsilon}}}\right]^{-n}}$$

$$\ln\left(\frac{\|z^n\|}{E_n(f)}\right)^{\frac{1}{n}}\geq \ln
\frac{\mu_1-\varepsilon}{\mu_2+\varepsilon}+\frac{1}{n}\ln
\left[{1-\frac{\mu_2+\varepsilon}{(\ln^{(q-2)}n)^{\frac{1}{\varepsilon}}}}\right]+\frac{1}{\varepsilon}
\ln ^{(q-1)}n$$

\begin{equation}\label{6}\frac{\ln\left(\frac{\|z^n\|}{E_n(f)}\right)}{n\ln^{(q-1)}n}\geq \frac{1}{\ln^{(q-1)}n}\cdot\ln
\frac{\mu_1-\varepsilon}{\mu_2+\varepsilon}+\frac{1}{n{\ln^{(q-1)}n}}\ln
\left[{1-\frac{\mu_2+\varepsilon}{(\ln^{(q-2)}n)^{\frac{1}{\varepsilon}}}}\right]+\frac{1}{\varepsilon}.
\end{equation} Отсюда

$$\liminf_{n\rightarrow\infty}\frac{\ln\left(\frac{\|z^n\|}{E_n(f)}\right)}{n\ln^{(q-1)}
n}\geq\frac{1}{\varepsilon}\,\,\text{ и, следовательно, }\,\,
\alpha_q=\limsup_{n\rightarrow\infty}\frac{n\ln^{(q-1)}
n}{\ln\frac{\|z^n\|}{E_n(f)}}\leq\varepsilon,$$ что противоречит
условию теоремы, значит $\rho_q>0$.

Выберем $\varepsilon \in (0,\mu_1)\cap (0,\rho_q).$ Из того, что
$$\rho_q=\limsup_{n\rightarrow\infty}\frac{n\ln^{(q-1)}
n}{-\ln|c_n|}$$ следует, что существует $ N_{\varepsilon}\in
\mathbb{N}$ такое, что при $ n
> N_{\varepsilon}$ выполняется неравенство   $
|c_n|<(\ln^{(q-2)}n)^{-\frac{n}{\varepsilon+\rho_q}}$. Будем считать
$N_{\varepsilon}$ столь большим, чтобы выполнялись неравенства:
$\|z^n\|\leq(\mu_2+\varepsilon)^n$,
$\,\,\|z^n\|\geq(\mu_1-\varepsilon)^n$ и
${\mu_2+\varepsilon}<{(\ln^{(q-2)}n)^{\frac1 {\varepsilon+\rho_q}}}$
при $n
> N_{\varepsilon}$. Тогда при $ n
> N_{\varepsilon}$

$$E_n(f)\leq \| \sum_{k=n}^\infty c_kz^k\| \leq
\sum_{k=n}^\infty |c_k|\cdot\|z^k\|\leq \sum_{k=n}^\infty
|c_k|\cdot(\mu_2+\varepsilon)^k\leq $$
\begin{equation}\label{6.1}\leq
\sum_{k=n}^\infty
\left[\frac{\mu_2+\varepsilon}{(\ln^{(q-2)}n)^{\frac{1}{\varepsilon+\rho_q}}}\right]^k
=
\frac{1}{1-\frac{\mu_2+\varepsilon}{(\ln^{(q-2)}n)^{\frac{1}{\varepsilon+\rho_q}}}}\cdot\left[\frac{\mu_2+\varepsilon}{(\ln^{(q-2)}n)^{\frac{1}{\varepsilon+\rho_q}}}\right]^n\,.\end{equation}
Следовательно,
$$\frac {\|z^n\|}{E_n(f)}\geq (\mu_1-\varepsilon)^n \cdot \left({1-\frac{\mu_2+\varepsilon}{(\ln^{(q-2)}n)^{\frac{1}{\varepsilon+\rho_q}}}}\right)\cdot\left[\frac{\mu_2+\varepsilon}{(\ln^{(q-2)}n)^{\frac{1}{\varepsilon+\rho_q}}}\right]^{-n}\,,$$

$$\ln\left(\frac{\|z^n\|}{E_n(f)}\right)^{\frac{1}{n}}\geq \ln
\frac{\mu_1-\varepsilon}{\mu_2+\varepsilon}+\frac{1}{n}\ln
\left[{1-\frac{\mu_2+\varepsilon}{(\ln^{(q-2)}n)^{\frac{1}{\varepsilon+\rho_q}}}}\right]+\frac{1}{\varepsilon+\rho_q}
\ln ^{(q-1)}n\,,$$

\begin{equation}\label{6}\frac{\ln\left(\frac{\|z^n\|}{E_n(f)}\right)}{n\ln^{(q-1)}n}\geq \frac{1}{\ln^{(q-1)}n}\cdot\ln
\frac{\mu_1-\varepsilon}{\mu_2+\varepsilon}+\frac{1}{n{\ln^{(q-1)}n}}\ln
\left[{1-\frac{\mu_2+\varepsilon}{(\ln^{(q-2)}n)^{\frac{1}{\varepsilon+\rho_q}}}}\right]+\frac{1}{\varepsilon+\rho_q}\,.
\end{equation} Отсюда

\begin{equation}\label{7}\rho_q+\varepsilon\geq \frac{n\ln^{(q-1)}
n}{\ln\frac{\|z^n\|}{E_n(f)}}\cdot
\left({1+\frac{\varepsilon+\rho_q}{\ln^{(q-1)} n}\cdot \ln
\frac{\mu_1-\varepsilon}{\mu_2+\varepsilon}+\frac{\varepsilon+\rho_q}{n\ln^{(q-1)}
n}\cdot}\ln
\left[{1-\frac{\mu_2+\varepsilon}{(\ln^{(q-2)}n)^{\frac{1}{\varepsilon+\rho_q}}}}\right]\right)\,.\end{equation}

Устремляя $n\rightarrow\infty$ получим, что
$\rho_q+\varepsilon\geq\alpha_q$. Ввиду произвольности выбора
$\varepsilon>0$ получаем $\rho_q\geq\alpha_q$ . Учитывая это и
неравенство (\ref{5}) имеем, что $\rho_q=\alpha_q$. Таким образом,
достаточность доказана.

\begin{flushleft}
    \emph{Необходимость.}
\end{flushleft}

Пусть  $f\in X$ целая функция конечного порядка $\rho_q$, т.е.

$$\rho_q=\limsup_{n\rightarrow\infty}\frac{n\ln^{(q-1)} n}{-|c_n|}\,.$$

Положим  $$\alpha=\limsup_{n\rightarrow\infty}\frac{n\ln^{(q-1)}
n}{\ln\frac{\|z^n\|}{E_n(f)}}\,.$$

Из леммы \ref{lemma1.3} следует, что $\alpha\geq\rho_q$. Рассуждая
как и при доказательстве достаточности, можем утверждать, что для
произвольного $\varepsilon\in (0,\mu_1)$ найдется $N_\varepsilon:
|c_n|<(\ln^{(q-2)}n)^{-\frac{n}{\varepsilon+\rho_q}}$,
$(\mu_1-\varepsilon)^n\leq\|z^n\|_D\leq(\mu_2+\varepsilon)^n$ и
 ${\mu_2+\varepsilon}<{(\ln^{(q-2)}n)^{\frac1 {\varepsilon+\rho_q}}}$ при $n > N_{\varepsilon}$.

Рассуждая как при доказательстве (\ref{6.1}) и (\ref{7}), получим
$$\rho_q+\varepsilon\geq \frac{n\ln^{(q-1)}
n}{\ln\frac{\|z^n\|}{E_n(f)}}\,\times$$
$$\times
\left({1+\frac{\varepsilon+\rho_q}{\ln^{(q-1)} n}\cdot \ln
\frac{\mu_1-\varepsilon}{\mu_2+\varepsilon}+\frac{\varepsilon+\rho_q}{n\ln^{(q-1)}
n}\cdot}\ln
\left[{1-\frac{\mu_2+\varepsilon}{(\ln^{(q-2)}n)^{\frac{1}{\varepsilon+\rho_q}}}}\right]\right)\,.$$

После предельного перехода при $n\rightarrow\infty$ в силу
произвольности выбора $\varepsilon>0$ получаем $\rho_q\geq\alpha$.
Учитывая обратное неравенство $\rho_q\leq\alpha$, имеем
$\rho_q=\alpha$. Таким образом, необходимость доказана.
\end{proof}

Докажем теорему \ref{th2}:

\begin{proof}

\begin{flushleft}
В случае $q=2$ наша теорема совпадает с теоремой \ref{th1.5},
доказанной в \cite{Dmz}.

    \emph{Достаточность.}
\end{flushleft}
Рассмотрим случай $q=3,4...\,$. Пусть $f\in X$ удовлетворяет условию
(\ref{4}) теоремы \ref{th2}, где $\rho_q$ и $\sigma_q$ некоторые
положительные числа. Тогда из (\ref{4}) следует справедливость
условия (\ref{3}) теоремы \ref{th1}, следовательно, $f$ - целая
функция порядка $\rho_q$. Пусть тип $f$ равен $\alpha$. Докажем, что
$\alpha=\sigma_q$. Из формулы для определения типа целой функции
через тейлоровские коэффициенты:

\begin{equation}\label{8}
\alpha=\limsup_{n\rightarrow\infty}\ln^{(q-2)}n\cdot|c_n|^\frac{\rho_q}{n},\quad
q=3,4,...\quad.
\end{equation}

С учетом леммы \ref{lemma1.3} имеем $\alpha\leq\sigma_q$. Докажем
обратное неравенство.

Из (\ref{8}) следует, что для произвольного $ \varepsilon>0$
существует $ N_\varepsilon \in \mathbb{N}$, для которого выполняются
неравенства $|c_n|<\left({\frac{\varepsilon+\alpha}{\ln
^{(q-2)}n}}\right)^{\frac{n}{\rho_q}}$ и
$\left({\frac{\varepsilon+\alpha}{\ln
^{(q-2)}n}}\right)^{\frac{1}{\rho_q}}\cdot(\mu+\varepsilon)<1$ при
всех $n>N_\varepsilon$ .

Оценим наилучшее приближение функции $f$ сверху:

$$E_n(f)\leq \left\| \sum_{k=n}^\infty c_kz^k\right\| \leq \sum_{k=n}^\infty
\left({\frac{\varepsilon+\alpha}{\ln
^{(q-2)}k}}\right)^{\frac{k}{\rho_q}}\cdot(\mu+\varepsilon)^k \leq$$

$$\leq\sum_{k=n}^\infty
\left(\left({\frac{\varepsilon+\alpha}{\ln
^{(q-2)}n}}\right)^{\frac{1}{\rho_q}}\cdot(\mu+\varepsilon)\right)^k\leq$$
\begin{equation}\label{15}
 \leq\left({1-\frac{C}{(\ln^{(q-2)}n)^\frac1
{\rho_q}}}\right)^{-1}\cdot \left({\frac{\varepsilon+\alpha}{\ln
^{(q-2)}n}}\right)^{\frac{n}{\rho_q}}\cdot(\mu+\varepsilon)^{n},
\end{equation}

 где
$C=(\mu+\varepsilon)(\varepsilon+\alpha)^{\frac{1}{\rho_q}}.\quad$
Из (\ref{15}) находим
$$\varepsilon+\alpha\geq \left(\frac{E_n(f)}{\|z^n\|}\right)^\frac{\rho_q}{n}\cdot\ln^{(q-2)}n\cdot\left({1-\frac{C}{(\ln^{(q-2)}n)^\frac1
{\rho_q}}}\right)^{\frac{\rho_q}{n}}\cdot\|z^n\|^{\frac{\rho_q}{n}}\cdot
(\mu+\varepsilon)^{-\rho_q} .$$

Последовательно устремив $n\rightarrow\infty$ и
$\varepsilon\rightarrow 0$, получим

$$\varepsilon+\alpha\geq \sigma_q\cdot\left(\frac{\mu}{\mu+\varepsilon}\right)^{\rho_q},\quad\alpha\geq\sigma_q,$$ что завершает доказательство достаточности.

\begin{flushleft}
    \emph{Необходимость.}
\end{flushleft}

Пусть $f\in X$ - целая функция обобщенного порядка $\rho_q, \quad
\rho_q\in(0;\infty)$. Обозначим ее обобщенный тип $\alpha$.
Аналогично доказательству достаточности, используя лемму
\ref{lemma1.3} и теорему \ref{th1.4}, можно показать справедливость
неравенства $\alpha\geq\sigma_q$. Чтобы доказать обратное
неравенство $\alpha\leq\sigma_q$ нужно повторить соответствующие
рассуждения из доказательства достаточности.

\end{proof}

Доказательство следствия 1 может быть получено применением соответствующих рассуждений на стр. 1132  работы  \cite{WacGir}.

Справедливость следствия 2 вытекает из того, что  при любом $r>0 \quad \rho_q(f)=\rho_q(f_r)$,  $\quad r^\frac{2}{p} \|f_r\| \leq \|f\|_\Omega \leq R^\frac{2}{p} \|f_R\| \,$  и теоремы  \ref{th2}.

Сведения об авторах:

1. Двейрин Михаил Захарович

 Донецкий национальный университет,

Кафедра математического анализа и дифференциальных уравнений

г. Донецк, 83055 Ул. Университетская, 24

Е-mail: matem47@mail.ru

тел. +38(062)-2972953, +380667075599

2. Левадная Антонина Сергеевна

 Донецкий национальный университет,

Кафедра математического анализа и дифференциальных уравнений

г. Дружковка,  Ул. Фестивальная, 4

 Е-mail: last.dris@mail.ru

тел. +380955268505

\end{document}